\documentstyle[aps,12pt]{revtex}
\begin{document}
\author{J. C. Stenerson}
\title{A Simple RevTeX Article}
\date{}

\mathstrut

\mathstrut

\quad \quad \quad \quad \quad {\bf \ Elementary Derivation of Spitzer's
Asymptotic Law for }

\quad \quad \quad \quad \quad {\bf Brownian Windings and Some of its
Physical Applications}

\mathstrut

\quad \quad \quad \quad \quad

\mathstrut

\quad \quad \quad \quad \quad \quad \quad \quad \quad \quad \quad \quad
\quad Arkady L. Kholodenko\footnote{%
375 H.L.Hunter Laboratories ,Clemson University,
\par
Clemson, SC 29634-1905 . e-mail :string@mail.clemson.edu}

\mathstrut

\mathstrut

\quad \quad \quad \quad \quad \quad \quad \quad \quad \quad \quad \quad
\quad \quad \quad \quad Abstract

\mathstrut Simple derivation of Spitzer's asymptotic law for Brownian
windings (Trans.Am.Math.Soc. {\bf 87} ,187 (1958)) is presented along with
its generalizations. These include the cases of planar Brownian walks
interacting with a single puncture and Brownian walks on a single truncated
cone with variable conical angle interacting with the truncated conical tip.
Such situations are typical in the theories of quantum Hall effect and 2+1
quantum gravity respectively. They also have some applications in polymer
physics. Extension of these results to the the multiple punctured case is
also briefly discussed. It is technically associated with some results known
in the context of string and conformal field theories and theories of
quantum chaos. \quad \quad \quad \quad \quad \quad \quad \quad \quad \quad
\quad

\mathstrut

\mathstrut


\mathstrut

\mathstrut

\newpage

In 1958 ,based on earlier work by Paul Levy [1] ,Spitzer [2] had obtained
the asymtotic probability distribution P(x) for the winding angle $\theta $
for the planar Brownian motion.If z(t)=x(t)+iy(t) is 2 -dimensional Wiener
process,then it is of interest to study distributions of $\left| z(t)\right| 
$ and $\arg z(t)$=$\theta (t)$ .For large times t , Spitzer obtained his
famous Cauchy-type distribution for $\theta (t)$ given by 
\begin{equation}
P(x=\frac{2\theta }{\ln t})dx=\frac{1}{\pi }\frac{1}{1+x^{2}}dx\text{ .} 
\eqnum{1}
\end{equation}
Eq.(1) is obtained under the assumption that the random walker begins his
travel at some point z in z-plane other than the origin , i.e. x=y=0 .
Then,the angle $\theta $ in measured with respect to the line which joins z
with the origin.This problem is of interest in polymer physics [3-5] since
it represents the benchmark problem for study of entanglements .In
mathematical literature the same phenomenon is described in terms of
recurrence and transience.For example,it is well known [6] ,that 1- and 2-
dimensional Brownian motions are recurrent (that is the random walk visits
time and again its starting point) while 3- dimensional motion is transient
(that is there is a nonzero probability that the walker will not return to
the origin).For the random walk on once- punctured z-plane, Lyons and McKean
[7] had demonstrated that the walk is recurrent while for the twice-
punctured plane the walk is transient. In the language of polymers this
means,that the polymer lying in z-plane will not be entangled with another
polymer, placed perpendicular to this plane , while it will become entangled
if there are at least 2 polymers which intersect z-plane at 2 distinct
points.The planarity of the above problem is actually not to essential as it
was explained in Ref.[5]. In the case of quantum mechanics the once puctured
plane problem is directly associated with the Aharonov-Bohm (A-B) effect
[8]. The A-B effect in the presence of 2 punctures was studied in
Ref.[9].The methods of Ref.[9] cannot be generalized to the case of more
than two punctures and do not provide any information about the recurrence
and/or transience.At the same time, the methods used by McKean and Lyons [7]
can be used for the case of more than 2 punctures but are not widely known
in physics literature. They had been recently mentioned in Ref.[5] in
connection with some topological problems arising in polymer physics.In
physics literature random walks on multiply punctured plane were extensively
studied in connection with problems related to the quantum Hall effect (QHE)
and anyonic superconductivity[8,10] while in mathematics literature the same
problem was recently extensively studied by Pitman and Yor [11,12].To our
knowlege,no attempt had been made to establish connections between these two
formalisms (see also Ref.[7]) . In this brief note we would like to make the
first step towards this comprehensive goal.

Let us begin with the well-known expression for the distribution function
for the planar random walk given by 
\begin{equation}
G({\bf r}_{1}\text{,{\bf r}}_{2};t)=\frac{1}{2\pi t}\exp \{-\frac{({\bf r}%
_{1}-{\bf r}_{2})^{2}}{2t}\}\text{ ,}  \eqnum{2}
\end{equation}
where {\bf r}$_{1}$({\bf r}$_{2})=\{x_{1}(x_{2}),y_{1}(y_{2})\}.$With
respect to the origin in z-plane the polar system of coordinates can be
used.In this system of coordinates Eq.(2) can be rewritten as 
\begin{equation}
G(r_{1},r_{2},\Delta \theta \text{ ; t)=}\frac{1}{2\pi t}\exp \{-\frac{%
r_{1}^{2}+r_{2}^{2}}{2t}\}\sum\limits_{m=-\infty }^{\infty }e^{im\Delta
\theta }I_{m}(z)\text{ ,}  \eqnum{3}
\end{equation}
,where $\Delta \theta =\theta _{1}^{{}}-\theta _{2},z=2r_{1}r_{2}/t$ and $%
I_{m}(z)=I_{-m}(z)$ is the modified Bessel's function.The above distribution
function can be used either for study of the radial or the angular
distributions or both.Suppose,we are interested in the angular distribution
function only (in view of Eq.(1) ).Then,using Eq.(3) ,it is convenient to
introduce the normalized distribution function defined according to the
following prescription:

\begin{equation}
f(z,\Delta \theta )=\frac{G(r_{1},r_{2},\Delta \theta \text{ ; t)}}{%
G(r_{1},r_{2},0\text{ ; t)}}=\frac{1}{I_{0}(z)}\sum\limits_{m=-\infty
}^{\infty }e^{im\Delta \theta }I_{\left| m\right| }(z)\text{ .}  \eqnum{4}
\end{equation}
The Fourier transform of such defined distribution function can now be
obtained in a standard way as 
\begin{equation}
f(z,\alpha )=\int\limits_{-\infty }^{\infty }d\Delta \theta \text{ e}%
^{-i\alpha \Delta \theta }f(z,\Delta \theta )=\frac{I_{\left| \alpha \right|
}(z)}{I_{0}(z)}\text{ .}  \eqnum{5}
\end{equation}
Let us now choose $r_{2}=\hat{r}\sqrt{t}+r_{1}.$ This choice is motivated by
known scaling properties of Brownian motion[13]. Then , for large t one
obtains :$z\simeq 2r_{1}\hat{r}/\sqrt{t}$ .For fixed $\hat{r}$ and $r_{1}$
and t$\rightarrow \infty $ one surely expects $z\rightarrow 0.$ This
observation allows us to use known asymptotic expansion for $I_{\left|
\alpha \right| }(z)$ for small $z^{\prime }s$ with the result for $%
f(z,\alpha )$ \quad ( valid for small $z^{\prime }s$ or large t's ): 
\begin{equation}
f(z,\alpha )\approx \exp \{-\frac{\left| \alpha \right| }{2}\ln t\}\text{ .}
\eqnum{6}
\end{equation}
The inverse Fourier transform of Eq.(6) leads us to the result given by
Eq.(1),i.e. $f(z,\Delta \theta )=P(x)$ where x=$2\Delta \theta /\ln t.$

Thus obtained result can be easily generalized now.For example,instead of
considering random walks on the flat once punctured plane we can consider
the same problem on the surface of a cone.This type of problem is of
interest in connection with the study of 2+1 quantum gravity [14,15]. It is
also known [15,16] ,that the above conical problem is equivalent to the
planar random walk problem in the wedge(the conical angle is simply related
to that of the wedge). In the most general case our walk may be allowed to
interact with the edges of the wedge[2] or,in the case of the wedge angle
equal to 2$\pi $ ,with the puncture located at the origin.The analogue of
the distribution function ,Eq.(3), is known to be [16,17]: 
\begin{equation}
G(r_{1},r_{2},\Delta \theta \text{ ; t)=}\frac{1}{\beta t}\exp \{-\frac{%
r_{1}^{2}+r_{2}^{2}}{2t}\}\sum\limits_{m=-\infty }^{\infty }e^{i2\pi
(m+\delta )\frac{\Delta \theta }{\beta }}I_{\frac{2\pi }{\beta }\left|
m+\delta \right| }(z)\text{ .}  \eqnum{7}
\end{equation}
For $\beta =2\pi $ and $\delta =0$ Eq.(7) reduces to Eq.(3) as required.The
wedge angle $\beta $ lies between 0 and 2$\pi $ while the statistics-
changing parameter $\delta $ is responsible for the polymer puncture
interactions as it is explained in Refs.[5,18] or for the interaction with
the flux tube if the magnetic language is being used [8].

By analogy with Eq.(4),we obtain, 
\begin{equation}
f_{\beta }^{\delta }(z,\Delta \theta )=\frac{1}{I_{\frac{2\pi }{\beta }%
\left| \delta \right| }(z)}\sum\limits_{m=-\infty }^{\infty }e^{i2\pi m\frac{%
\Delta \theta }{\beta }}I_{\frac{2\pi }{\beta }\left| m+\delta \right| }(z)%
\text{ .}  \eqnum{8}
\end{equation}
Upon Fourier transforming this expression we obtain, 
\begin{equation}
f_{\beta }^{\delta }(z,\alpha )=\frac{I_{\left| \alpha +\frac{2\pi }{\beta }%
\delta \right| }(z)}{I_{\frac{2\pi }{\beta }\left| \delta \right| }(z)}\text{
\quad .}  \eqnum{9}
\end{equation}
Repeating the same chain of arguments which had lead us to Eq.(6) we obtain
now, 
\begin{equation}
f_{\beta }^{\delta }(z,\alpha )\simeq \exp \{-\frac{1}{2}(\left| \alpha +%
\frac{2\pi }{\beta }\delta \right| -\frac{2\pi }{\beta }\left| \delta
\right| )\ln t\}\text{ .}  \eqnum{10}
\end{equation}
To perform the inverse Fourier transform of Eq.(10) is nontrivial.Indeed,we
have 
\begin{equation}
f_{\beta }^{\delta }(\Delta \theta ,t)=\frac{1}{2\pi }\int\limits_{-\infty
}^{\infty }d\alpha \exp \{i\alpha \Delta \theta \text{ -}\frac{1}{2}(\sqrt{%
\alpha ^{2}+\left( \frac{2\pi }{\beta }\delta \right) ^{2}}-\frac{2\pi }{%
\beta }\left| \delta \right| )\ln t\}\text{ .}  \eqnum{11}
\end{equation}
The integrals of this type are known in the context of quantum field theory
[19] and had been also used recently in polymer physics problems[20].By
introducing new variable $\alpha =\left| \frac{2\pi \delta }{\beta }\right|
\sinh \varphi $ into Eq.(11) it is transformed into 
\begin{equation}
f_{\beta }^{\delta }(\Delta \theta ,t)=\frac{\delta }{\beta }\exp \{\frac{%
\pi }{\beta }\left| \delta \right| \ln t\}\int\limits_{-\infty }^{\infty
}d\varphi \cosh \varphi \exp \{-\frac{1}{2}[\cosh \varphi ]\frac{2\pi }{%
\beta }\left| \delta \right| \ln t+i\Delta \theta \frac{2\pi }{\beta }\left|
\delta \right| \sinh \varphi \}\text{ .}  \eqnum{12}
\end{equation}
The exponent inside of the integral in Eq.(12) can be transformed as
follows: 
\[
-\frac{1}{2}[\cosh \varphi ]\frac{2\pi }{\beta }\left| \delta \right| \ln
t+i\Delta \theta \frac{2\pi }{\beta }\left| \delta \right| \sinh \varphi =-%
\sqrt{a^{2}+\omega ^{2}}\cosh (\varphi +\varphi _{0})\text{ ,} 
\]
where $a=\frac{\pi }{\beta }\left| \delta \right| \ln t$ and $\omega =\Delta
\theta \frac{2\pi }{\beta }\left| \delta \right| $ so that $\cosh \varphi
_{0}=\frac{a}{\sqrt{a^{2}+\omega ^{2}}}$ and $\sinh \varphi _{0}=\frac{%
-i\omega }{\sqrt{a^{2}+\omega ^{2}}}$ .Use of these results in Eq.(12)
allows us to rewrite it in the equivalent form , 
\begin{eqnarray}
f_{\beta }^{\delta }(\Delta \theta ,t) &=&\frac{\delta \frac{2\pi }{\beta }%
\left| \delta \right| \ln t}{\beta \sqrt{a^{2}+\omega ^{2}}}\exp \{\frac{\pi 
}{\beta }\left| \delta \right| \ln t\}\int\limits_{0}^{\infty }d\varphi
\cosh \varphi \exp [-\sqrt{a^{2}+\omega ^{2}}\cosh \varphi ]  \nonumber \\
&\equiv &\Phi _{\beta }^{\delta }(t)K_{1}(\sqrt{a^{2}+\omega ^{2}})\text{ ,}
\eqnum{13}
\end{eqnarray}
where $K_{1}(x)$ is the modified Bessel's function with known asymptotic
expansions: $K_{1}(x)\simeq \frac{1}{x}$ for $x\rightarrow 0$ and $%
K_{1}(x)\simeq \sqrt{\frac{\pi }{2x}}e^{-x}$ for $x\rightarrow \infty .$
Using these expansions ,the following asymptotic results for the
distribution function $f_{\beta }^{\delta }(\Delta \theta ,t)$ are obtained :

1.$\delta \rightarrow 0$ and $\beta $ is fixed and nonzero.In this case we
recover Spitzer's law , Eq.(1) as required.

2.$\delta \rightarrow \infty $ but $\Delta \theta $ is finite.Then,we obtain
(x=$\frac{2\Delta \theta }{\ln t})$ , 
\begin{equation}
f_{\beta }^{\delta }(x)dx\simeq \frac{1}{2}\sqrt{\frac{\delta }{2\beta }\ln t%
}\frac{dx}{(1+x^{2})^{\frac{3}{4}}}\exp \{-\left| \frac{\beta \ln t}{8\pi
\delta }\right| x^{2}\}\text{ .}  \eqnum{14}
\end{equation}
Obtained result is in complete accord (up to numerical prefactors) with that
obtained in Ref.[21] using methods of conformal field theory and,more
recently, in Ref.[22] ,by direct numerical simulations.

3.$\delta $ and $\beta $ are finite but $\Delta \theta \rightarrow \infty $
and $\Delta \theta \gg \ln t$ .Then, we obtain , 
\begin{equation}
f_{\beta }^{\delta }(x)dx\simeq \frac{1}{2}\sqrt{\frac{\delta }{2\beta }\ln t%
}\exp \{\left| \frac{\pi \delta }{\beta }\right| \ln t\}\exp \{-\left|
x\right| \left| \frac{\pi \delta }{\beta }\right| \ln t\}\frac{dx}{\left|
x\right| ^{\frac{3}{2}}}\text{ .}  \eqnum{15}
\end{equation}
This result is in formal agreement with that obtained by Rudnick and Hu [4]
and Desbois [23].In Ref.[4] only the leading exponential factor was obtained
while the result,Eq.(15), differs from that obtained in Ref.[23] by an extra
time-dependent factor .No connections with conical singularities or QHE
problems were made in either Ref.[4] or Ref.[23].

With the results just obtained ,we are ready now to make connections with
works of Ito and McKean [24] and Lyons and McKean[7] (see also Ref.[25]).We
begin by stating the result proved by Levy[1] and refined by others [13].
Let z(t) be some planar Brownian motion which started ,say,at z =0 (at t =
0). Then , the motion obtained with help of some analytic function $f(z)$ is
also Brownian. This is equivalent to saying that the planar brownian motion
is conformally invariant. Let us illustrate this fact in the example of once
punctured plane {\bf R}$^{2}-{\bf 0.}$

The diffusion equation (in dimensionless units) on {\bf R}$^{2}-{\bf 0}$ can
be written in a usual form as 
\begin{equation}
\frac{\partial f}{\partial t}=\frac{1}{4}\left( \frac{\partial ^{2}}{%
\partial x^{2}}+\frac{\partial ^{2}}{\partial y^{2}}\right) f\text{ ,} 
\eqnum{16}
\end{equation}
where the factor $\frac{1}{4}$ is chosen for convenience.It is useful to
rewrite Eq.(16) in terms of complex variables.Simple calculation produces 
\begin{equation}
\frac{\partial f}{\partial t}=\frac{\partial ^{2}f}{\partial z\partial \bar{z%
}}\text{ .}  \eqnum{17}
\end{equation}
Since this equation is not defined for z = \={z} = 0,we introduce new
complex variable $w$ through z=$\exp w$ .Unlike z, our new variable $w=u+iv$
is defined for the entire complex $w-$plane.In terms of $w$-variable Eq.(17)
can be rewritten as 
\begin{equation}
\frac{\partial f}{\partial t}=\frac{1}{2}e^{-2u}\left( \frac{\partial ^{2}}{%
\partial u^{2}}+\frac{\partial ^{2}}{\partial v^{2}}\right) f  \eqnum{18}
\end{equation}
Obtained equation describes Brownian motion on a simply- connected covering
space $\tilde{M}$ which is Riemann surface for the logarithmic function.The
result just obtained coincides with that discussed in the book by Ito and
McKean,e.g.see page 280 of ref.[24].Ito and McKean argue (without explicit
demonstration) that ,based on the results of P.Levy ,the diffusion Eq.(18)
can be converted into standard form given by Eq.(16) by replacing Brownian
time t by another (actually random) time T which is properly chosen. Let us
demonstrate how this can actually be done.To this purpose, let us consider
the Langevin-type equation written in the form of Ito : 
\begin{equation}
dx(s)=a(x(s))ds+\sigma (x(s))dw(s)\text{ .}  \eqnum{19}
\end{equation}
The corresponding backward Kolmogorov-Fokker-Planck equation can be written
now as 
\begin{equation}
\frac{\partial }{\partial s_{0}}P(x,s\mid x_{0},s_{0})=-a\frac{\partial }{%
\partial x_{0}}P+\frac{1}{2}\sigma ^{2}\frac{\partial ^{2}}{\partial
x_{0}^{2}}P\text{ .}  \eqnum{20}
\end{equation}
Introduce ''new'' time T(s) according to equation 
\begin{equation}
s=\int\limits_{0}^{T(s)}dtg(x(t))\text{ .}  \eqnum{21}
\end{equation}
Therefore, $ds=g(x(T))dT.$Substitution of this result into Eq.(19) and use
of Ito stochastic calculus [13] produces, 
\begin{equation}
dx(T)=a(x(T))g(x(T))dT+\sqrt{g(x(T))}\sigma (x(T))dw(T)\text{ .}  \eqnum{22}
\end{equation}
New diffusion coefficient can be selected now ( in view of Eq.(20))as 
\begin{equation}
D_{new}=\frac{1}{2}g\sigma ^{2}\text{ .}  \eqnum{23}
\end{equation}
Looking at Eq.(21) and selecting

\mathstrut

\begin{equation}
s=\int\limits_{0}^{T(s)}dte^{2u(t)}  \eqnum{24}
\end{equation}
produces , in view of Eq.(23) , $D_{new}=\frac{1}{2}$ . Eq.(24) is in
complete agreement with the result of Ito and McKean [21] where it was
presented without derivation . Surely, T(s) is random time \quad so that
solution of Eq.(16) will depend upon random time.This is somewhat
inconvenient.To correct this inconvenience it is useful to think about
calculable observables.In view of Eq.(5) ,we are interested in the
Fourier-transformed angular distribution function $f(z,\alpha )$ .It can be
shown [24],that this task is equivalent to finding averages of the type 
\begin{equation}
f(z,\alpha )=\left\langle e^{i\alpha \Delta \theta }\right\rangle
=\left\langle \exp \{-\frac{\alpha ^{2}}{2}\int\limits_{0}^{T}d\tau g(x(\tau
))\}\right\rangle ,  \eqnum{25}
\end{equation}
where $\left\langle ...\right\rangle $ denotes the averaging with help of
Gaussian-like propagator for the {\bf free }''particle''. Eq.(25) is a
special case of famous Feynman-Kac formula[26] . The transition from first
to second average in Eq.(25) is associated with the
Cameron-Martin-Girsanov-type formula for changes of variables inside of path
integrals[13,27].Eq.(25) is associated with Schrodinger-like differential
equation.In view of Eq.s(21) and (24) ,it can be easily demonstrated ,using
standard methods of quantum mechanics [28] , that such obtained
Shrodinger-like equation is just that for the modified Bessel functions $%
,e.g.$ $I_{m}$ (x),etc .Its solution , indeed, leads us to Eq.(5).

Based on these results ,immediate generalizations are possible.For
instance,if we would be looking for the distribution function $f(z,\alpha )$
for 2- (or more)punctured plane,we would need the uniformizing function
which will help us to write the diffusion equation on the simply connected
universal covering surface $\tilde{M}$ .This surface is known to be the
Riemann surface for the punctured torus [25].Classical motion on such
surface is chaotic and quantum description of such motion was considered in
Ref.[29] where no attempts were made to relate it to Spitzer's
results.Already for two punctures study of transience and recurrence is
highly nontrivial [7,25]. For two or more punctures the task of finding the
uniformizing function can be reduced to that of finding all solutions of the
corresponding Fuchsian-type equations[30] which are identical in the form to
that known in string and conformal field theories for the correlation
functions[31].It remains a challenging problem to extend Spitzer-like
results to multiply connected surfaces,to connect these results with those
of Pitman and Yor[11,12] and with those known in the theories of QHE [8,10].

{\bf Note added in proof}. When this work was completed, two additional
recent results came to our attention. In Ref.[32] reader can find more up to
date(as compared to Pitman an Yor's) references on Brownian windings while
Ref.[33] contains detailed applications of the above ideas to QHE

\newpage 

\quad \quad \quad \quad \quad \quad \quad \quad \quad \quad

\quad \quad \quad \quad \quad \quad \quad \quad \quad \quad \quad \quad
\quad {\bf References}

\mathstrut

[1] P.Levy , {\it Processus Stochastiques et Movement Brownien }(Paris ,1948
)

[2] F.Spitzer ,Trans .Am.Math.Soc. {\bf 87 },187 (1958 )

[3] S.Edwards ,Proc.Phys.Soc. {\bf 91 }, 513 (1967)

[4] J.Rudnick and Y.Hu , Phys.Rev. Lett. {\bf 60 },712 (1988); and also
J.Phys.A{\bf 20,}

\quad 4421 (1987)

[5] A.Kholodenko and Th.Vilgis , {\it Some Geometrical and Topological
Problems in}

\quad {\it Polymer} {\it Physics} , Physics Reports {\bf 298}(5-6) (1998).

[6] F.Spitzer ,{\it Principles of Random Walk (}Van Nostrand,
Princeton,N.J.,1964)

[7] T.Lyons and H.McKean ,Adv.in Math.{\bf 51 }212 (1984)

[8] A.Shapere and F.Wilczek ,{\it Geometric Phases in Physics} (World

\quad Scientific,Singapore,1990)

[9] S.-Wu Qian, Zhi-Y. Gu and Guo-Q.Xie ,J.Phys.A {\bf 30,}1273 (1997{\bf )}

[10] F.Wilczek , {\it Fractional Statistics and Anyon Superconductivity}
(World

\quad \quad Scientific,Singapore ,1990)

[11] J.Pitman and M.Yor ,Ann Prob. {\bf 14 },733 91986)

[12] J.Pitman and M.Yor ,Ann Prob. {\bf 17 },965 (1989)

[13] R.Durrett, {\it Brownian Motion and Martingales} (Wadsworth Inc.
Belmont,CA ,1984)

[14] S.Deser ,R.Jackiw and G.'t Hooft ,Ann Phys. {\bf 152 },220 (1994)

[15] A.Kholodenko , J.Math. Phys. (1998) ,submitted

[16] J.Dowker ,J.Phys.A {\bf 10 },115 (1977)

[17] S.Deser and R.Jackiw ,Comm.Math.Phys.{\bf \ 118 },495 (1988)

[18] A.Kholodenko ,J.Douglas and D.Bearden,Phys.Rev.E{\bf 49 },2206 (1994)

[19] N.Bogoliubov and D.Shirkov ,{\it Introduction to the Theory of
Quantized Fields}

\quad \quad (Wiley ,N.Y.,1980)

[20] A.Kholodenko and Th.Vilgis ,Phys.Rev.E{\bf 50, }1257(1994)

[21] B.Duplantier and H.Saleur ,Phys.rev.Lett. {\bf 60 }2343 (1988)

[22] T.Prellberg and DB.Drossel, Phys.Rev.E {\bf 57 },2045 (1998)

[23] J.Desbois , J.Phys.A {\bf 23 },3099 (1990)

[25] H.McKean and d.Sullivan ,adv.math. {\bf 51 },203 (1984)

[26] M.Freidlin and A.Wentzell, {\it Random Perturbations of Dynamical
Systems}

\quad \quad (Springer-Verlag,Berlin ,1984)

[27] I.Gelfand and A.Yaglom ,J.Math.Phys.{\bf 1},48 (1961)

[28] R.Feynman and A.Hibbs ,{\it Quantum Mechanics and Path Integrals(}%
McGraw Hill{\it ,}

{\it \ \quad \quad }N.Y.,1965{\it )}

[29] T.Shigehara, Phys.rev.E {\bf 50 },4357 (1994)

[30] J.Hempel , Bull.London Math.Soc. {\bf 20 },97 (1988)

[31] P.Gisparg and G.Moore in {\it Recent Directions in Particle Theory ,}%
pages 277-469

\quad \quad (World Scientific,Singapore,1993)

[32] Z.Shi , Ann.Probability {\bf 26 },112 (1998).

[33] A.Carey ,K.Hannabuss ,V.Mathai and P.McCann ,Comm.Math.Phys. 

\quad \quad {\bf 190,}629(1998).

\end{document}